\documentclass{article}
\usepackage{nips10submit_e,times}
\usepackage{algorithm, algorithmic}
\usepackage{tensorStyle}
\usepackage{subfig}
\usepackage{graphicx}
\usepackage{epstopdf}
\usepackage{url}
\usepackage{paralist}

\title{Making Tensor Factorizations Robust to Non-Gaussian Noise\thanks{This paper was a contributed presentation at the NIPS Workshop on Tensors, Kernels, and Machine Learning, Whistler, BC, Canada, December 10, 2010.}}

\author{
Eric C.~Chi\footnotemark[1] \\
Department of Statistics\\
Rice University\\
Houston, TX 15213 \\
\texttt{echi@rice.edu} \\
\And
Tamara G.~Kolda\footnotemark[2]\\
Sandia National Laboratories \\
Livermore, CA 94551-9159 \\
\texttt{tgkolda@sandia.gov} \\
}

\nipsfinalcopy 

\begin{document}

\maketitle

\renewcommand{\thefootnote}{\fnsymbol{footnote}}
\footnotetext[1]{This work was funded by DOE grant DE-FG02-97ER25308.}
\footnotetext[2]{This work was funded by the applied mathematics program at the U.S. Department of Energy. Sandia National Laboratories is a multi-program laboratory managed and operated by Sandia Corporation, a wholly owned subsidiary of Lockheed Martin Corporation, for the U.S. Department of Energy's National Nuclear Security Administration under contract DE-AC04-94AL85000.} 
\renewcommand{\thefootnote}{\arabic{footnote}}

\begin{abstract}
Tensors are multi-way arrays, and the CANDECOMP/PARAFAC (CP) tensor factorization has found application in many different domains. The CP model is typically fit using a least squares objective function, which is a maximum likelihood estimate under the assumption of i.i.d.\@ Gaussian noise. We demonstrate that this loss function can actually be highly sensitive to non-Gaussian noise. Therefore, we propose a loss function based on the 1-norm because it can accommodate both Gaussian and grossly non-Gaussian perturbations. We also present an alternating majorization-minimization algorithm for fitting a CP model using our proposed loss function.
\end{abstract}

\section{Introduction}
The CANDECOMP/PARAFAC (CP) tensor factorization can be considered a higher-order generalization of the matrix singular value decomposition \cite{Carroll1970, Harshman1970} and has many applications. 
The canonical fit function for the CP tensor factorization is based on the Frobenius norm, meaning that it is a maximum likelihood estimate (MLE) under the assumption of additive i.i.d.\@ Gaussian perturbations. It turns out, however, that this loss function can be very sensitive to violations in the Gaussian assumption. However, many other types of noise are relevant for CP models. For example, in fMRI neuroimaging studies, movement by the subject can lead to sparse high-intensity changes that are easily confused with brain activity \cite{Friston1996}. Likewise, in foreground/background separation problems in video surveillance, a subject walking across the field of view represents another instance of a sparse high intensity change \cite{Li2004}. In both examples, there is a relatively large perturbation in magnitude that affects only a relatively small fraction of data points; we call this artifact noise. These scenarios are particularly challenging because the perturbed values are on the same scale as normal values (i.e., true brain activity signals and background pixel intensities). Consequently, there is a need to explore factorization methods that are robust against violations in the Gaussian assumption. In this paper, we consider a loss based on the 1-norm which is known to be robust or insensitive to gross non-Gaussian perturbations \cite{Huber2009}.

Vorobyov et al.\@ previously described two ways of solving the least 1-norm CP factorization problem based on a linear programming and weighted median filtering  \cite{Vorobyov2005}. Our method differs in that we use a majorization-minimization (MM) strategy  \cite{Hunter2004}. Like \cite{Vorobyov2005} our method performs block minimization. An advantage of our approach is that each block minimization can be split up into many small and independent optimization problems which may scale more favorably with a tensor's size.

Throughout, we use the following definitions and conventions. All vectors are column vectors. The transpose of the $i^{\rm th}$ row of a matrix $\M{A}$ is denoted by $\Mr{A}{i}$.
The {\em order} of a tensor is the number of dimensions, also known as ways or modes. {\em Fibers} are the higher-order analogue of matrix rows and columns. A fiber is defined by fixing every index but one. A matrix column is a mode-$1$ fiber and a matrix row is a mode-$2$ fiber. The mode-$n$ matricization of a tensor $\T{X} \in \mathbb{R}^{I_1 \times I_2 \times \cdots \times I_N}$ is denoted by $\Mz{X}{n}$ and arranges the mode-$n$ fibers to be the columns of the resulting matrix.

The rest of this paper is organized as follows. The robust iterative algorithm is derived in Section~\ref{sec:MM}. In Section~\ref{sec:simulation} we compare CPAL1 and the standard CP factorizations by alternating least squares (CPALS) in the presence of non-Gaussian perturbations on simulated data. Concluding remarks are given in Section~\ref{sec:conclusion}.

\section{Majorization-minimization for tensor factorization}
\label{sec:MM}
MM algorithms have been applied to factorization problems previously \cite{Lee2001, Bro2002, Dhillon2006}.  The idea is to convert a hard optimization problem (e.g., non-convex, non-differentiable) into a series of simpler ones (e.g., smooth convex), which are easier to minimize than the original. To do so, we use majorization functions, i.e., $h$ majorizes $g$ at $\V{x}\in \mathbb{R}^n$ if $h(\V{u}) \geq g(\V{u})$ for all $\V{u} \in \mathbb{R}^n$ and $h(\V{x}) = g(\V{x})$.

Given a procedure for constructing a majorization, we can define the MM algorithm to find a minimizer of a function $g$  as follows. Let $\Vn{x}{k}$ denotes the $k^{\rm th}$ iterate.
\begin{inparaenum}[(1)]
  \item Find a majorization $h(\cdot | \Vn{x}{k})$ of $g$ at $\Vn{x}{k}$.
  \item Set $\Vn{x}{k+1} = \arg \min_{\V{x}} h(\V{x} | \Vn{x}{k})$.
  \item Repeat until convergence.
\end{inparaenum}
This algorithm always takes non-increasing steps with respect to $g$. Moreover, sufficient conditions for the MM algorithm to converge to a stationary point are well known \cite{Lange2010}. Specifically, the MM iterates will converge to a stationary point of $g$ if $g$ is continuously differentiable, coercive in the sense that all its level sets must be compact, and all its stationary points are isolated; and
 the function $h(\V{x} | \V{y})$ is jointly twice continuously differentiable in $(\V{x}, \V{y})$ and is strictly convex in $\V{X}$ with $\V{y}$ fixed.

\subsection{Solving the $\ell_1$ regression problem by an MM algorithm}
\label{sec:LADr}

We now derive an appropriate majorization function for approximate $\ell_1$ regression; this is subsequently used for our robust tensor factorization.  The unregularized majorization can be found in \cite{Lange2010}.
Given a vector $\V{Y} \in \mathbb{R}^I$ and a matrix $\M{M} \in \mathbb{R}^{I \times J}$, we search for a vector $\V{u} \in \mathbb{R}^J$  that minimizes the loss $L(\V{u}) = \sum_i \lvert \VE{r}{i}(\V{u}) \rvert$ where 
$\VE{r}{i}(\V{u}) = \VE{Y}{i} - \Mr{M}{i}\Tra\V{u}$.
Note that $L(\V{u})$ is not smooth and may not be strictly convex if $\M{M}$ is not full rank. Therefore, we instead consider the following smoothed and regularized version to $L(\V{u})$:
\begin{equation}
\label{eq:el1min}
L_{\epsilon, \mu}(\V{u}) = \sum_{i=1}^I \sqrt{\VE{r}{i}(\V{u})^2 + \epsilon} + \frac{\mu}{2}\lVert \V{u} \rVert^2,
\end{equation}
where $\epsilon$ and $\mu$ are small positive numbers.
In this case, $L_{\epsilon,\mu}(\V{u})$ at $\Vtilde{u} \in\mathbb{R}^J$ is majorized by
\begin{equation}
\label{eq:least_el1_maj}
h_{\epsilon,\mu}(\V{u} | \Vtilde{u}) =  \sum_{i=1}^I \left \{ \sqrt{\VE{r}{i}(\Vtilde{u})^2 + \epsilon} + \frac{\VE{r}{i}(\V{u})^2- \VE{r}{i}(\Vtilde{u})^2}{2\sqrt{\VE{r}{i}(\Vtilde{u})^2 + \epsilon}} \right \} +  \frac{\mu}{2}\lVert \V{u} \rVert^2. \\
\end{equation}
Both the loss $L_{\epsilon,\mu}$ and its majorization $h_{\epsilon,\mu}$ meet the sufficient conditions that guarantee convergence of the MM algorithm to
a stationary point of $L_{\epsilon,\mu}$. Since $L_{\epsilon, \mu}$ is strictly convex and coercive, it has exactly one stationary point. Thus, the MM algorithm converges to the global minimum of $L_{\epsilon,\mu}$.
After some simplification, the iterate mapping can be expressed as 
\begin{equation}
\label{eq:least_el1_MMb}
\Vn{u}{m+1} = \arg\underset{\V{u}}{\min}   \left\{\sum_{i=1}^I \frac{\VE{r}{i}(\V{u})^2}{\sqrt{\VE{r}{i}(\Vn{u}{m})^2+ \epsilon}} +  \frac{\mu}{2}\lVert \V{u} \rVert^2 \right \}.
\end{equation}
Let $\Mn{W}{m} \in\mathbb{R}^{I \times I}$ be the diagonal matrix with
 $(\Mn{W}{m})_{ii} =  (r_i(\Vn{u}{m})^2+ \epsilon)^{-1/2}$. Then the minimization problem (\ref{eq:least_el1_MMb}) is a regularized weighted least squares problem with a unique solution, i.e.,
\begin{equation}
\label{eq:least_el1_MM_sol}
\Vn{u}{m+1} = \arg\underset{\V{u}}{\min} \; \left \{ (\V{Y} - \M{M}\V{u})\Tra \Mn{W}{m} (\V{Y} - \M{M}\V{u}) +  \frac{\mu}{2}\lVert \V{u} \rVert^2 \right \} 
= (\M{M}\Tra \Mn{W}{m} \M{M} + \mu \bm{\mathbf{I}})^{-1} \M{M} \Tra\Mn{W}{m} \V{y}.
\end{equation}
\subsection{Tensor factorization using $\ell_1$ regression}
We now derive CPAL1 for a 3-way tensor $\T{X}$ of size $I_1 \times I_2 \times I_3$ (it is straightforward to generalize the algorithm to tensors of arbitrary size). To perform a rank-$R$ factorization
we minimize $L_{\epsilon,\mu}(\M{A},\M{B},\M{C})$ which is the regularized approximate 1-norm of the difference between $\T{X}$ and its rank-$R$ approximation,
where $\M{A} \in \mathbb{R}^{I_1\times R}, \M{B} \in \mathbb{R}^{I_2\times R},  \M{C} \in \mathbb{R}^{I_3 \times R}$:
\begin{equation*}
L_{\epsilon,\mu}(\M{A},\M{B},\M{C}) =\sum_{i_1,i_2,i_3} \sqrt{\left (\TE{X}{i_1i_2i_3} - \sum_{r=1}^R \ME{A}{i_1r}\ME{B}{i_2r}\ME{C}{i_3r} \right )^2 + \epsilon} +  \frac{\mu}{2}\left (\lVert
\M{A} \rVert_F^2 
+ \lVert
\M{B} \rVert_F^2
+ \lVert
\M{C} \rVert_F^2
\right).
\end{equation*}
In a round-robin fashion, we repeatedly update one factor matrix while holding the other two fixed. Note that the mode-1 matricization of
the rank-$R$ approximation is $\M{A}(\M{C} \Khat \M{B})\Tra$ where $\Khat$ denotes the Khatri-Rao product \cite{Kolda2009}. Then the subproblem of updating $\M{A}$ for a fixed $\M{B}$ and $\M{C}$ is
\begin{equation}
\label{eq:sep}
\min_{\M{A}}\sum_{i=1}^{I_1}  \sum_{j=1}^{I_2 I_3} \sqrt{\left ( \left(\Mz{X}{1} - \M{A}(\M{C} \Khat \M{B})\Tra\right)_{ij} \right)^2 + \epsilon} +  \frac{\mu}{2}\lVert \M{A} \rVert_F^2.
\end{equation}
This minimization problem is separable in the rows of $\M{A}$, and the optimization problem for a given row is an $\ell_1$ regression problem. Thus, we can apply the update rule (\ref{eq:least_el1_MM_sol})
with $\V{y}$ equal to the $i^{\rm th}$ row of $\Mz{X}{1}$ and $\M{M} = \M{C} \Khat \M{B}$. The other two subproblems are solved analogously.

\section{Simulation Experiment \label{sec:simulation}}
We compare the results of CPAL1 with CPALS implemented in the Tensor Toolbox \cite{Bader2010} in the presence of Gaussian and artifact noise. We created $3$-way tensors, $\T{X}' \in \Real^{50\times 50 \times 50}$ of rank-$5$ as follows. We first generated random factor matrices $\M{A}, \M{B}, \M{C} \in \Real^{50 \times 5}$ where the matrix elements were the absolute values of i.i.d.\@ draws from a standard Gaussian. The $ijk^{\rm th}$ entry of the noise free tensor $\T{X}$ was then set to be $\sum_{r=1}^R \ME{A}{i_1r}\ME{B}{i_2r}\ME{C}{i_3r}$. Then to each $\T{X}$ we added dense Gaussian noise and artifact outliers. All random variables we describe were independently drawn. We generated an artifact tensor $\T{P}$ as follows.  A fraction $\eta$ of the tensor entries was selected randomly. We then assigned to each of the selected entries a value drawn from a Gamma distribution with shape parameter $50$ and scale parameter $1/50$. All other entries were set to $0$.
For the dense Gaussian noise tensor $\T{Q}$, the entries $\TE{Q}{ijk}$ were i.i.d. draws from a standard Gaussian. 
The tensor $\T{X}'$ was obtained by adding the noise and artifact tensors to $\T{X}:$
\begin{equation*}
\T{X}' = \T{X} + \gamma \frac{{\lVert \T{X}\rVert_F}}{\lVert \T{P}\rVert_F} \T{P}  + 0.1 \frac{{\lVert \T{X}\rVert_F}}{\lVert \T{Q}\rVert_F} \T{Q} 
\end{equation*}
for $\eta = 0.1, 0.2$ and $\gamma = 0.5, 1.0, 1.5$ and $2.0$. For all combinations of $\eta$ and $\gamma$ the scaled values of $\TE{Q}{ijk}$ were less than the largest value of $\T{X}$.

For every pair $(\eta, \gamma)$ we performed $100$ rank-$5$ factorizations under the two methods. For CPAL1 computations we set $\epsilon = 10^{-10}$ and $\mu = 10^{-8}$. Initial points for all tests were generated using the $n$-mode singular vectors of the tensor (i.e., the \verb nvecs  command in the Tensor Toolbox).
To assess the goodness of a computed factorization we calculated the factor match score (FMS) between the estimated and true factors \cite{Acar}.
The FMS ranges between 0 and 1; an FMS of $1$ corresponds to a perfect recovery of the original factors.

\begin{figure}[t]
\centering
\subfloat[The FMS distribution under different combinations of $\eta$ and $\gamma$.]
{\label{fig:boxplots}
\includegraphics[scale=0.575]{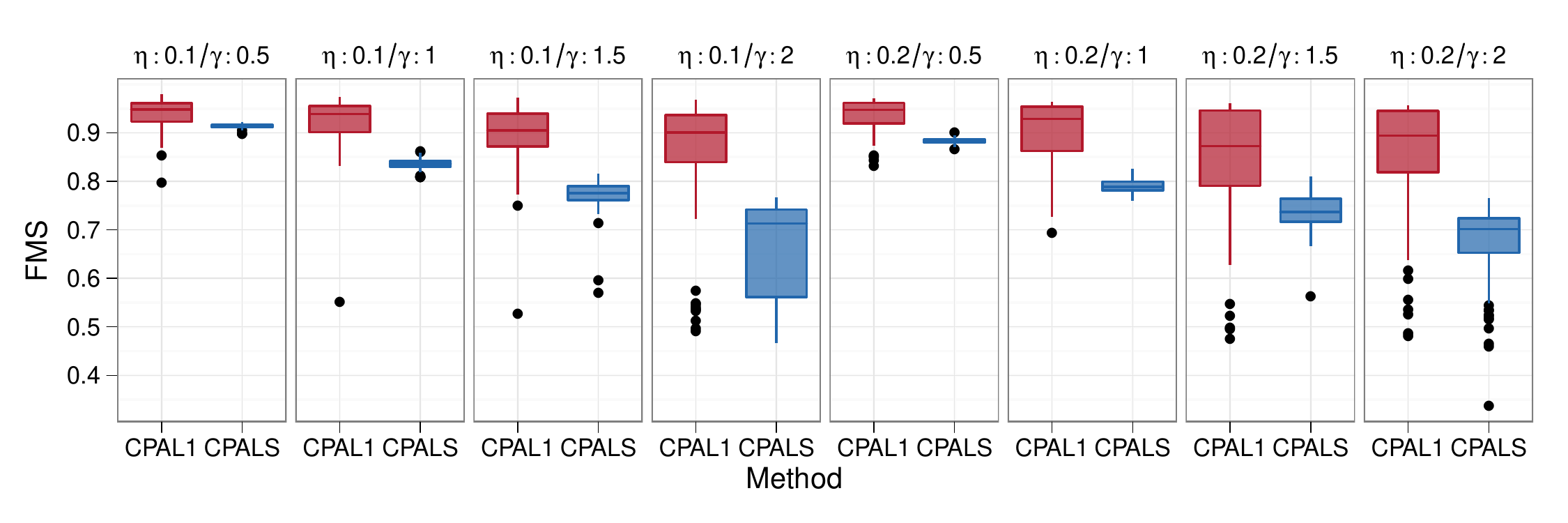}
} \\
\subfloat[A comparison of a single recovered factor column for a replicate when $\eta=0.2$ and $\gamma=2$. Here the FMS was $0.91$ and $0.64$ for CPAL$1$ and CPALS respectively. Factor columns were normalized for comparison.]
{\label{fig:factors}
\includegraphics[scale=0.575]{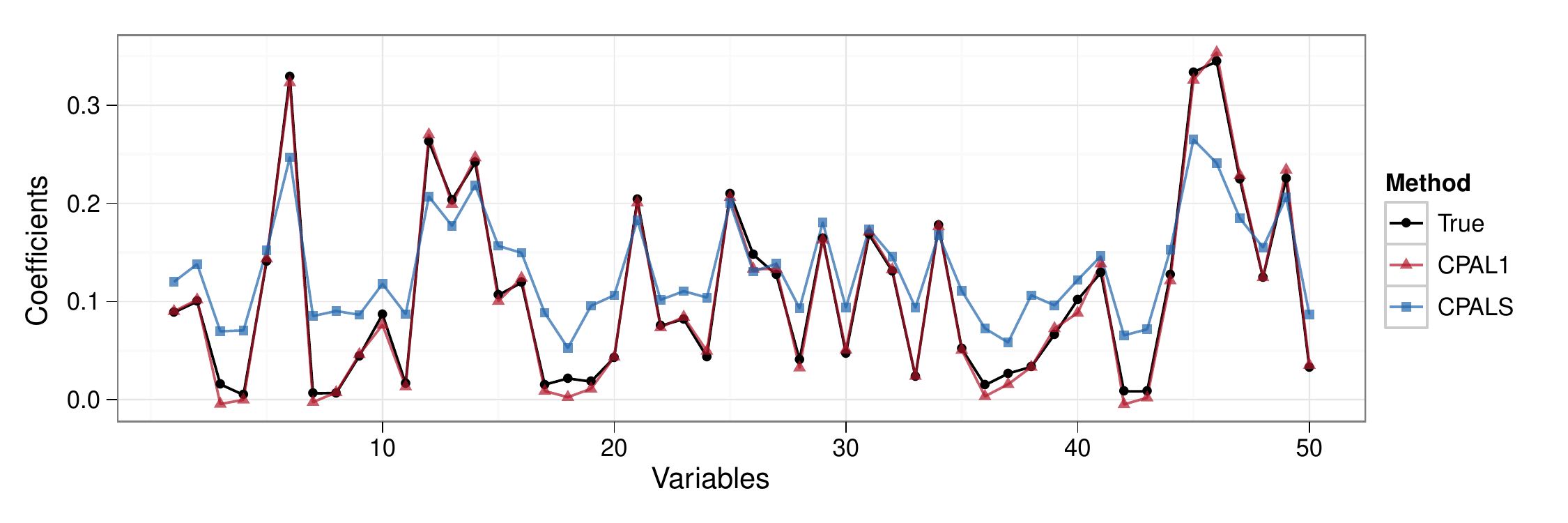}
}
\caption{Panel~\ref{fig:boxplots} shows on average that CPALS factorizations are sensitive to artifact noise. Panel~\ref{fig:factors} provides a close up comparison 
of the differences between the two methods for single recovered column in an instance where CPALS is less accurate.}
\end{figure}

Figure~\ref{fig:boxplots} shows boxplots of the FMS under both methods. The scores for CPALS decreased as the contribution of non-Gaussian noise increased. In contrast regardless of the noise distributions applied CPAL1 tended to recover the true factorization with the exception of occasionally finding local minima, 

Figure~\ref{fig:factors} compares one column of one recovered factor matrix when $\eta = 0.2$ and $\gamma=2$ for the two methods. In this instance
the CPALS factorization has some trouble recovering the true factor column.  In this example the FMS was $0.91$ and $0.64$ for CPAL$1$ and CPALS respectively. The median CPALS FMS was about $0.7$, so the example shown is somewhat typical.  The factorization is not terrible qualitatively, but the errors in the Factor 2 estimates do fail to capture details that CPAL1 solution does.

\section{Conclusion \label{sec:conclusion}}
We derived a robust tensor factorization algorithm based on an approximate 1-norm loss. In comparisons with methods using an LP solver we found that our method performed slightly faster on tensors of similar size to those factored in the simulation experiments of this paper (not shown). We suspect the performance gap may widen depending on the size of the tensor. Indeed, to factor an arbitrary tensor of size $I_1\times \cdots \times I_N$
the LP update for the $i^{\rm th}$ factor matrix would be an optimization problem over $2 \prod_{n \not= i} I_n + R$ parameters. In contrast, the $i^{\rm th}$ factor matrix update consists of $I_i$ independent $\ell_1$ minimizations over $R$ parameters. Moreover, the independence of these minimizations present speed-up opportunities through parallelization.

Our simulations demonstrated that there are non-Gaussian noise scenarios in which the quality of CPALS solutions suffer while those of CPAL1 tend to be insensitive to the presence of non-Gaussian noise. In simulation studies not shown we have seen that not all non-Gaussian perturbations cause noticeable degradation in the CPALS factorization. Conversely, there are situations when CPAL1 struggles as much as CPALS in the presence of artifact noise, e.g. when the data tensor is sparse as well. We conjecture that CPAL1 is most suited to handle artifact noise when the data tensor is dense. Finding an alternative to the 1-norm loss for sparse data with non-Gaussian noise is a direction for future research.

\bibliographystyle{siam}

\end{document}